\documentclass[letterpaper,10pt,oneside,onecolumn]{amsart}
\usepackage{amsmath, amssymb}
\usepackage[left=1in,top=1in,right=1in,bottom=1in]{geometry}

\newcommand{\EE}{\mathbb E}
\newcommand{\F}{\mathcal F}

\newcommand{\PP}{\mathbb P}

\newcommand{\R}{\mathbb R}
\newcommand{\Z}{\mathbb Z}

\newcommand{\Beuc}{B^{\mathrm{E}}}

\newcommand{\smin}{\setminus}

\newcommand{\SPD}{\operatorname{SPD}}

\newcommand{\E}{\mathrm{e}}				% Euler's constant.
\newcommand{\sD}{\, \mathrm{d}}				% Differential operator with space.
\newcommand{\oo}{\infty}

\numberwithin{equation}{section}
\theoremstyle{definition}
\newtheorem{env_thm}{Theorem}[section]
\newtheorem{env_lem}[env_thm]{Lemma}
\newtheorem{env_cor}[env_thm]{Corollary}
\newtheorem{env_pro}[env_thm]{Proposition}
\newtheorem{env_exa}[env_thm]{Example}

\title{A Shape Theorem for Riemannian First-Passage Percolation}
\author{T. LaGatta and J. Wehr}
\begin{document}
	
	\begin{abstract}
 		Riemannian first-passage percolation (FPP) is a continuum model, with a distance function arising from a random Riemannian metric in $\R^d$.  Our main result is a shape theorem for this model, which says that large balls under this metric converge to a deterministic shape under rescaling.  As a consequence, we show that smooth random Riemannian metrics are geodesically complete with probability one.
	\end{abstract}

	\maketitle

	\section{Introduction}
	
	\subsection{Motivation}
	
	We introduce Riemannian first-passage percolation (FPP) as a model of random geometry in the continuum.  Our work is motivated by standard FPP, where one builds a random distance function on $\Z^d$ from independent, identically-distributed ``passage times'' over bonds of the lattice.  We instead construct our distance function on $\R^d$ from a random Riemannian metric.  The main result of the paper is a shape theorem for Riemannian FPP:  if $B_t$ is the random Riemannian ball of radius $t$, then $\tfrac{1}{t}B_t$ tends toward a limiting shape.  From this it follows that a random Riemannian metric is geodesically complete with probability one.
	
	Hammersley and Welsh \cite{hammersley1965fpp} introduced standard FPP in 1965 in order to model fluid flow through porous media.  Consider the $d$-dimensional lattice $\Z^d$ with $d \ge 2$.  Let $\{t_b\}$ be a family of independent, identically distributed, non-negative random variables, indexed by bonds (nearest-neighbor edges) $b$ of the lattice.  For any $z,z' \in \Z^d$, define
		$$d(z,z') = \inf_{\gamma} \sum_{b \in \gamma} t_b,$$
	where the infimum is taken over all lattice paths $\gamma$ connecting $z$ to $z'$.  This $d$ is a random distance function on $\Z^d$.  For a very good introduction to standard FPP, see Howard \cite{howard2004mfp}.  
	
	Consider the random distance $d(0, n\E_1)$ between the origin and the point $n\E_1 = (n,0,\dots,0)$.  One wishes to study the asymptotic behavior of this quantity as $n \to \infty$.  In \cite{kingman1968ets}, Kingman formulated his famous subadditive ergodic theorem in order to prove the basic result of FPP:  provided the passage times have finite mean, there exists a non-random constant $\mu_{\E_1}$, such that 
		$$\lim_{n\to\oo} \tfrac{1}{n} d(0,n\E_1) = \mu_{\E_1}$$
	almost surely and in $L^1$.  The same is clearly true for all coordinate axes and, more generally, for each direction $v \in S^{d-1}$, there exists a non-random constant $\mu_v$, such that
		\begin{equation} \label{kingmandiscrete}
			\lim_{n\to \oo} \tfrac{1}{n}d(0,\widetilde{nv}) = \mu_v \end{equation}
	almost surely and in $L^1$, where $\widetilde{nv} \in \Z^d$ is the nearest lattice point to $nv$.  The constant $\mu_v$ is non-zero provided that the probability that $t_b = 0$ is less than the critical percolation probability for $\Z^d$ \cite{howard2004mfp}.

	The shape theorem of Cox and Durrett \cite{cox1981slt} is a stronger result.  Consider
		$$\tilde B_t = \{z \in \Z^d : d(0,z) \le t\},$$
	the random ball of radius $t$ in $\Z^d$.  This is a lattice object, so we ``inflate'' it to get a continuum one:  for $z \in \Z^d$, let $C_z = [z-1/2,z+1/2)^d$ be the unit cube centered at $z$ in $\R^d$, and let 
		$$B_t = \bigcup_{z \in \tilde B_t} C_z.$$
	We define the rescaling $\tfrac{1}{t} B_t$ as the set of all points $x \in \R^d$ such that $tx \in B_t$.  The shape theorem says that there exists a non-random convex, compact set $A$, depending only on the distribution of $t_b$, such that $\tfrac{1}{t} B_t \to A$:  for all $\epsilon > 0$, with probability one, there exists a time $T$ such that if $t \ge T$, then
		$$(1-\epsilon)A \subseteq \tfrac{1}{t} B_t \subseteq (1+\epsilon)A.$$
	
	Just as standard FPP is a model of random geometry on the discrete lattice $\Z^d$, Riemannian FPP is a model of random geometry in the continuum.  We consider a random Riemannian metric $g$ on $\R^d$ whose distribution is translation-invariant, has finite-range dependence and satisfies certain moment conditions.  By the standard construction in Riemannian geometry, this defines a random distance function $d(x,y)$ in $\R^d$.  Kingman's theorem can again be applied to prove that for each unit vector $v$, there exists a non-random constant $\mu_v \ge 0$ such that $\tfrac{1}{n} d(0,nv) \to \mu_v$ a.s. and in $L^1$.  From the positive-definiteness of $g$, we prove Theorem \ref{posconst}:  $\mu_v > 0$ for all $v \in S^{d-1}$.  Consider the set
		$$A = \{x \in \R^d : |x| \le \mu_{x/|x|}^{-1} \},$$
	and the Riemannian ball of radius $t$ centered at the origin
		$$B_t = \{ x \in \R^d : d(0,x) \le t \}.$$
	Theorem \ref{shapethm} is the shape theorem for Riemannian FPP, which states that for all $\epsilon > 0$, with probability one, there exists a random time $T > 0$ such that if $t \ge T$, then
		$$(1-\epsilon)A \subseteq \tfrac{1}{t} B_t \subseteq (1+\epsilon)A.$$
	Consequently, we call $A$ the limiting shape of the model.  
	
	Theorem \ref{completeness} follows from the shape theorem:  with probability one, no curve $\gamma$, parametrized by Riemannian length, reaches infinity in finite time.  When the metric is further assumed to be smooth with probability one, then this is geometrically significant:  by the Hopf-Rinow theorem \cite{lee1997rmi}, this is equivalent to geodesic completeness of the metric.

	To prove our results, we need some technical estimates on the distance function $d$, which we obtain in Section \ref{disclemmas}.  To prove these, we discretize the continuum model:  to each point $z \in \Z^d$, we associate a certain value $X_z$ based on the Riemannian metric $g(x)$ over the unit cube $C_z = [z-1/2,z+1/2)^d$.  By treating $X_z$ as defining a dependent FPP model on the lattice, we prove some estimates of the ``entropy-energy type'' on $X_z$, and from these derive the desired estimates on $d$.

	\subsection{Geometry Background and Notation} 
	
	Before introducing any probabilistic structure, we introduce some geometric notation.  Consider $\R^d$ with $d \ge 2$ and the standard Euclidean coordinates.  Write
		$$\SPD = \{ \mbox{symmetric, positive-definite $d \times d$ real matrices} \},$$
	and let $g \in C(\R^d,\SPD)$ be a continuous matrix-valued function on $\R^d$ with values in $\SPD$.  $g$ defines a Riemannian structure on $\R^d$:  for tangent vectors $v, v' \in T_x \R^d$, we consider the inner product $\langle v, g(x) v' \rangle$.  For a single vector $v$, we denote by $\|v\| = \sqrt{\langle v, g(x) v \rangle}$ and $|v| = \sqrt{\langle v, v \rangle}$ the Riemannian and Euclidean lengths of $v$, respectively.  For a $C^1$-curve $\gamma : [a,b] \to \R^d$, we define the Riemannian and Euclidean lengths of $\gamma$ by
		$$R(\gamma) = \int_a^b \| \dot \gamma(t) \| \sD t \qquad \mathrm{and} \qquad L(\gamma) = \int_a^b | \dot \gamma(t) | \sD t,$$
	respectively.  We say that a curve is finite if it has finite Euclidean length; for our model, Theorem \ref{completeness} will imply that finite curves have finite Riemannian length.  The Riemannian distance between two points $x$ and $y$ is defined by
		$$d(x,y) = \inf_\gamma R(\gamma),$$
	where the infimum is over all $C^1$-curves $\gamma$ connecting $x$ to $y$.
	
	For a Riemannian metric $g$, we define the real, positive functions
		$$\Lambda(x) = \mbox{maximum eigenvalue of $g(x)$} \qquad \mathrm{and} \qquad \lambda(x) = \mbox{minimum eigenvalue of $g(x)$}.$$
	For any $K \subseteq \R^d$, define
		$$\Lambda(K) = \sup_{x\in K} \Lambda(x) \qquad \mathrm{and} \qquad \lambda(K) = \inf_{x\in K} \lambda(x).$$
	By the continuity and positivity of $g$, if $K$ is bounded then
		$$0 < \lambda(K) \le \Lambda(K) < \oo.$$
	For $z \in \Z^d$, let $C_z = [z - 1/2, z + 1/2)^d$ be the unit cube centered at $z$.  Write 
		$$\Lambda_z = \Lambda(C_z) \qquad \mathrm{and} \qquad \lambda_z = \lambda(C_z).$$

	\subsection{Riemannian FPP}

	Let $\Omega = C(\R^d, \SPD)$ and let $\F$ be the $\sigma$-algebra generated by cylinder sets.  Let $\PP$ be a translation-invariant probability measure on $(\Omega, \F)$ which has finite-range dependence, and consider a random Riemannian metric $g \in \Omega$ with distribution $\PP$.  The finite-range dependence means there exists some $R > 0$ such that if $|x-y|\ge R$, then $g(x)$ and $g(y)$ are independent.  Furthermore, suppose that $\Lambda_0$ has a finite moment-generating function.  That is,
		\begin{equation} \label{finmom}
			M(r) = \EE[ \E^{r \Lambda_0} ] < \oo \qquad \mbox{for all $r \in \R$,} \end{equation}
	where $\EE$ denotes expectation with respect to $\PP$.  By translation invariance, the family $\{\Lambda_z\}$ is identically distributed, and we refer to its generic element as $\Lambda$; similarly for $\{\lambda_z\}$ and $\lambda$.  By Chebyshev's inequality \cite{durrett1996probability}, \eqref{finmom} implies that $\Lambda$ and $\lambda$ have exponential tail decay:
		$$\PP( \lambda > u ) \le \PP( \Lambda > u ) \le M(r) \E^{-r u},$$
	for all $u > 0$ and $r > 0$.	
	
	We provide a concrete example.
	
	\begin{env_exa}
		Let $c : [0,\oo) \to \R$ be a compactly-supported covariance function (see \cite{gneiting2002csc} for examples).  Let $\xi : \R^d \to \R$ be a mean-zero, stationary, istropic Gaussian field with covariance function $c$; that is,
			$$\EE[\xi(x)] = 0 \qquad \mathrm{and} \qquad \EE[\xi(x)\xi(y)] = c(|x-y|)$$ 
		for all $x, y \in \R^d$.  The covariance $c$ must satisfy certain necessary and sufficient conditions \cite{talagrand1987regularity} for the field $\xi$ to be everywhere continuous with probability one; suppose this is the case.
		
		Let $g : \R^d \to \R$ be the diagonal matrix with entries
			$$g_{ii}(x) = \log(1 + \E^{\xi(x)}),$$
		for $1 \le i \le d$.  This is continuous and positive, so it suffices to show that the assumption \eqref{finmom} is satisfied.  
		
		\begin{env_pro}
			Assumption \eqref{finmom} is satisfied for this choice of $g$, and the shape theorem (Theorem \ref{shapethm}) applies.  Let $A$ be the limiting shape, defined in \eqref{limshape}.  The measure $\PP$ is isotropic, so Corollary \ref{isotropic} implies that $A$ is a Euclidean ball.  Furthermore, if the field $\xi$ is $C^1$ with probability one, then Corollary \ref{realized} implies that the metric $g$ is geodesically complete.
		\end{env_pro}
		\begin{proof}

		Assume for simplicity that $c(0) = 1$.  The Gaussian concentration inequality \cite{adler07} implies that
			$$\PP(\Lambda_0 > u) = \PP(\sup \xi > \log(e^u - 1) ) \le \exp(-\log(\E^u-1)^2/2) \le 2 \E^{-u^2/2}.$$
		The fundamental theorem of calculus and Fubini's theorem imply that
			$$\EE \E^{r \Lambda_0} = \EE \! \left( 1 + \int_0^{\Lambda_0} r \E^{ru} \sD u \right) = 1 + \int_0^{\oo} r \E^{ru} \PP( \Lambda_0 > u ) \sD u \le 1 + 2 \int_0^{\oo} r \E^{ru} \E^{-u^2/2} \sD u,$$
		which is finite for all $r$.  Thus \eqref{finmom} is satisfied, and the results of this paper apply to the random Riemannian metric $g$.
		\end{proof}
	\end{env_exa}

	The random variables $\lambda_z$ and $\Lambda_z$ give rise to a dependent FPP model on sites of the lattice $\Z^d$.  In Section \ref{depfpp}, we prove some general estimates for dependent FPP, then in Section \ref{contapps} we apply these to estimates on our distance function $d$.  Our techniques are based on the energy-entropy methods of mathematical physics, where one shows that an event occurs with extremely low probability over one particular connected set (``high energy''), but sums this over all possible connected sets at the origin (``high entropy'').  One adjusts parameters in the problem so that this sum converges, then applies the Borel-Cantelli lemma.  A large-deviations estimate like \ref{finmom} is critical: the number of connected sets at the origin grows exponentially in $n$, the size of the sets, so the probabilities must decay exponentially in $n$ for the arguments to hold.

	In the standard FPP setting of passage times $t_b$ across bonds $b$, Cox and Durrett \cite{cox1981slt} prove that a necessary and sufficient condition for a shape theorem is that $\EE \min\{t_1, \dots, t_{2d}\}^d < \oo$, where $t_i$ are $2d$ independent copies of $t_b$.  Thus we believe that our assumption \eqref{finmom} is not the most general, and can be replaced by a finite moment estimate on $\Lambda$ instead to prove a more general result.  
	
	In probability theory, Kingman's subadditive ergodic theorem \cite{durrett1996probability} is used to prove that stationary, subadditive sequences obey laws of large numbers.  If $X_{n,m}$ is a non-negative, stationary sequence which satisfies $X_{n,m} \le X_{n,r} + X_{r,m}$, then $\tfrac{1}{n} X_{0,n}$ converges almost surely and in $L^1$.  Furthermore, if the sequence is ergodic, this convergence is to a non-random constant.  In our context, the sequence in question is $X_{n,m} = d(nv,mv)$ for a fixed unit vector $v$.  The subadditivity condition is exactly the triangle inequality for $d$, so Kingman's theorem implies that for each $v \in S^{d-1}$, there exists a non-random constant $\mu_v \ge 0$ such that
		\begin{equation} \label{kingman}
			\lim_{t\to\oo} \tfrac{1}{t} d(0,tv) = \mu_v \end{equation}
	almost surely and in $L^1$.  The constants $\mu_v$ may depend on the direction $v$, though if the measure $\PP$ is isotropic (rotationally-invariant) then $\mu_v = \mu$ will be independent of $v$.  We show in Theorem \ref{posconst} that $\mu_v > 0$ for all $v$.
	
	\begin{env_pro} \label{contmu}
	 	$\mu_v$ is a continuous function of $v$.
	\end{env_pro}
	\begin{proof}
		This remarkably short proof is due to Kesten \cite{kesten1180arp} (see his Proof of Theorem 1.7 on page 158).  First we show that the function $\mu_v$ is bounded above as a function of $S^{d-1}$.  Write $v = \sum_1^d v^i \E_i$, for the standard basis vectors $\E_i$ in $\R^d$.  We use the triangle inequality to bound $d(0,tv)$ by the sum of the distances between successive points $0$, $t v^1 \E_1$, $t v^1 \E_1 + t v^2 \E_2$ and so on until $tv$.  Translation-invariance implies
	 		$$\EE d(0,tv) \le \EE d(0,tv^1 \E_1) + \dots + \EE d(0,tv^d \E_d).$$
	 	Dividing by $t$ and taking the limit $t \to \oo$ gives
 			$$\mu_v = \lim_{t\to\oo} \tfrac{1}{t} \EE d(0,tv) \le \sum_{i=1}^d \lim_{t\to\oo} \tfrac{1}{t} \EE d(0, tv^i \E_d).$$
 		The terms which equal zero we may ignore; for the non-zero terms, we make the substitution $t' = |v^i| t,$ so that the right-hand side equals
 			$$\sum_{i=1}^d \lim_{t'\to\oo} \tfrac{|v^i|}{t'} \EE d(0, \pm t' \E_i) = \sum_{i=1}^d |v^i| \mu_{\E_i} \le d \max\{\mu_{\E_i}\},$$
 		as desired.  
 		 	 		
	 	Now, consider two different unit vectors $v$ and $v'$, and write $u = \tfrac{v-v'}{|v-v'|}$.  By the same arguments,
	 		$$|\mu_v - \mu_{v'}| \le \lim_{t\to\oo} \tfrac{1}{t} \EE d(tv,tv') = \mu_u |v - v'|.$$
	 	This tends to zero as $v' \to v$ since $\mu_u$ is bounded above.
	\end{proof}

	\section{Discretization Lemmas} \label{disclemmas}
	
	\subsection{Dependent FPP on a Lattice} \label{depfpp}

		For a continuous curve $\gamma$, we would like to introduce a discrete analogue $\Gamma$ on the lattice.  For example, $z \in \Gamma$ if $\gamma$ meets the cube $C_z$.  However, this set $\Gamma$ is not connected on $\Z^d$; consider the straight line from $0$ to $(1,1)$ in $\R^2$.  We get around this by modifying the familiar graph structure of $\Z^d$ to introduce a new lattice, which we call the $*$-lattice.  In this section, we prove some estimates for dependent FPP on the $*$-lattice, then in Section \ref{contapps} we apply these estimates to the continuum model.
		
		For $z \in \Z^d$, we write $z = (z^1, \dots, z^d)$.  We say that $z, z' \in \Z^d$ are $*$-adjacent if $\max_{1\le i\le d} (z - z')^i \le 1$.  The $*$-lattice is the graph with vertex set $\Z^d$, and edge set given by $*$-adjacency; that is, the usual lattice $\Z^d$ along with all the diagonal edges.
		
		We say that a set $\Gamma \subseteq \Z^d$ is $*$-connected if for all $z, z' \in \Gamma$, there is a path from $z$ to $z'$ along the $*$-lattice which remains in the set $\Gamma$.  Technically, that there is a finite sequence of $*$-adjacent points beginning with $z$ and ending with $z'$, all contained in $\Gamma$.  
		
		Let $S_n$ be the number of $*$-connected sets which contain the origin.
		
		\begin{env_lem} \label{numsets}
			There exists $\sigma$ such that $S_n \le \sigma^n$.  Obviously, $\sigma > 1$.  
		\end{env_lem}
		\begin{proof}
			Clearly, $\log S_n$ is a non-negative subadditive sequence:
				$$\log S_{n+m} \le \log S_n + \log S_m.$$
			By Fekete's lemma \cite{fekete1923verteilung}, there exists a constant $a \ge 0$ such that $\log S_n \le a n$ for all $n$.  Defining $\sigma = \E^a$ proves the result.
		\end{proof}
		
		Let $X_z$ be a stationary, non-negative random field on the $*$-lattice with finite-range dependence, and with a finite moment-generating function
			\begin{equation} \label{finmom2}
				M(r) = \EE[\E^{r X}] < \oo \qquad \mbox{for all $r \in \R$.} \end{equation}
		The finite-range dependence means that there is an integer $R \ge 1$ such that if $|z - z'| \ge R$, then $X_z$ and $X_{z'}$ are independent.
		
		If $\Gamma \subseteq \Z^d$ is a collection of lattice points, we write
			$$X(\Gamma) = \sum_{z \in \Gamma} X_z,$$
		and call this the passage time of $\Gamma$.  %If a path $\tilde\gamma$ meets a point $z$ multiple times, then the field $X_z$ is only counted once in the passage time calculuation.
		
		The following two lemmas can be thought of as spatial laws of large numbers.  The first says that for sufficiently large $n$, if there is a uniform bound on $X(\Gamma)/n$, then there is also a uniform bound on $|\Gamma|/n$.  The second lemma reverses the implication, though with different constants.
				
		\begin{env_lem} \label{passstep}
			Suppose $X_z$ additionally satisfies
				\begin{equation} \label{zeroatom}
					\PP(X = 0) < \sigma^{-(2R+1)^d}. \end{equation}
			For any $A > 0$ there is a non-random $B > 0$ such that, with probability one, for any sequence $a_n \in \Z^d$, there exists $N > 0$ such that for all $n \ge N$, if $\Gamma$ is a $*$-connected set which contains the point $a_n$ and $X(\Gamma) \le An$, then $|\Gamma| \le Bn.$
		\end{env_lem}
	
		\begin{env_lem} \label{upbound}
			For any $B > 0$ there is a non-random $C > 0$ such that, with probability one, for any sequence $a_n \in \Z^d$, there exists $N > 0$ such that for all $n \ge N$, if $\Gamma$ is a $*$-connected set which contains the point $a_n$ and $|\Gamma| \le Bn$, then $X(\Gamma) \le Cn.$
		\end{env_lem}
		
		Some assumption like \eqref{zeroatom} is necessary for Lemma \ref{passstep}.  Let $p = \PP(X = 0)$, and suppose that $p > p_c$, the critical probability for site percolation on the $*$-lattice.  By percolation theory, with probability one, the set
			$$\Gamma = \{ z \in \Z^d : X_z = 0 \}$$
		contains an infinite $*$-connected component.  That is, $|\Gamma| = \oo$ but $X(\Gamma) = 0$.  No such assumption is necessary for Lemma \ref{upbound}.  
		
		In this paper, our assumptions of non-negativity and \eqref{zeroatom} are stronger than necessary:  we apply these lemmas only to the positive fields $\lambda_z$ and $\Lambda_z$.  However, we anticipate these lemmas to be of independent use in future work on Riemannian FPP, where one may consider fields $X_z$ which take the value $0$ on a cube $C_z$ with small but non-zero probability.  For example, if $E_z$ is the event that $\Lambda_z \le h$ for a sufficiently large value of $h$, one may apply these lemmas to $X_z = 1_{E_z}$, the indicator function of $E_z$.  
		
		The generality of the sequence $a_n$ is needed for Lemma \ref{unifT}.  In the proof of that result, we fix a point $x \in \R^d$, and define the sequence $a_n = \widetilde{nx} \in \Z^d$ to be the nearest lattice point to $nx$.

		\begin{proof}[Proof of Lemma \ref{passstep}]
		In what follows, we assume that $\Gamma$ is a $*$-connected set.  Fix $A > 0$ and the sequence $a_n$.  Consider the events
			$$E_n = \{ \exists ~ \Gamma \mathrm{~such~that~} a_n \in \Gamma, ~ X(\Gamma) \le An,\mathrm{~and~} |\Gamma| > Bn \}.$$
		We claim that we can choose an integer $B$ large enough so that $\PP(E_n)$ decays exponentially in $n$.  From the Borel-Cantelli lemma it will follow that, with probability one, only finitely many of the events $E_n$ occur, which will prove the result.

		Let $z_i$ be a predetermined enumeration of $\Z^d$; for example, a spiral path beginning at $a_n$.  For any $*$-connected set $\Gamma$, we define $\Gamma' \subseteq \Gamma$ by proceeding along the sequence $z_i$ and including each point of $\Gamma$ at a distance at least $R$ away from the previous points chosen.  In the form of an algorithm:
			\begin{itemize}
				\item Let $i_1$ be the first index for which $z_{i_1} \in \Gamma$.  Let $g_1 = z_{i_1}$.
				\item Given $\{g_1, \dots, g_{j-1}\}$, let $i_j$ be the first index for which $z_{i_j} \in \Gamma$ and so that $|z_{i_j} - g_{j'}| > R$ for $1 \le j' < j$.  Let $g_j = z_{i_j}$.
			\end{itemize}
		Let $\Gamma' = \{g_1, g_2, \dots\}$.  If $\Gamma$ is finite, then so is $\Gamma'$.  Let 
			$$B(z,R) = \{z' \in \Z^d : |z-z'| \le R\}$$
		be the Euclidean ball of radius $R$ centered at $z$ in $Z^d$.  Note that $|B(z,R)| \le (2R+1)^d$.  To avoid this cumbersome factor $(2R+1)^d$ which appears frequently, for the remainder of this proof we write
			$$K = (2R+1)^d.$$
		By construction, the set $\Gamma$ is covered by taking balls around every point in $\Gamma'$:
			$$\Gamma \subseteq \bigcup_{z \in \Gamma'} B(z,R).$$

		Let $\Gamma$ be a $*$-connected set described in the event $E_n$, so that
			$$Bn < |\Gamma| \le \sum_{z \in \Gamma'} |B(z,R)| \le |\Gamma'| K.$$
		Thus, $|\Gamma'| > Bn/K$.
		
		We may assume that $\Gamma$ consists only of the first $Bn$ points it meets of the sequence $z_i$; the non-negativity of $X_z$ implies that the passage time $X(\Gamma)$ is still at most $An$.  Similarly, we assume $|\Gamma'| = \lfloor Bn/K \rfloor$, the integer part of $Bn/K$.  Furthermore, since $\Gamma' \subseteq \Gamma$,
			$$X(\Gamma') \le X(\Gamma) \le An.$$
		Thus the probability $\PP(E_n)$ is bounded by
			\begin{equation} \label{splitestimate}
				\PP \left( \exists ~\Gamma \mathrm{~s.t.~} a_n \in \Gamma, ~ X(\Gamma') \le An, \mathrm{~and~} |\Gamma| = Bn \right) \le \sum_{\Gamma} \PP \left( X(\Gamma') \le An \right), \end{equation}
		where the outer sum is taken over all $*$-connected sets $\Gamma$ for which $a_n \in \Gamma$ and $|\Gamma| = Bn$.  Lemma \ref{numsets} implies that the number of such sets is bounded by $\sigma^{Bn}$.
		
		The family of random variables $\{X_z\}_{\Gamma'}$ is independent, since the points $z \in \Gamma'$ are separated by distances at least $R$.  Let $\{X_i\}$ be $\lceil Bn/K \rceil$ independent copies of $X$.  The exponential Chebyshev inequality \cite{durrett1988lnp} implies that the right-hand side of \eqref{splitestimate} is bounded above by
			$$\sigma^{Bn} \, \PP \left( \sum_{i=1}^{\lceil Bn/K \rceil} X_i \le An \right) \le \sigma^{Bn} \, \E^{r An} \, \EE \left( \E^{-r \sum X_i} \right)$$
		for any $r > 0$.  Again by independence, if we write $M(-r) = \EE \E^{-r X}$, this is equal to
			$$\sigma^{Bn} \E^{r An} \left(\EE \E^{-r X} \right)^{\lceil Bn/K \rceil} \le \sigma^{Bn} \E^{r An} M(-r)^{Bn/K} = \left( \sigma \, \E^{r A / B} \, M(-r)^{1/K} \right)^{Bn}.$$
	
		By the bounded convergence theorem, as $r$ tends to infinity, $M(-r) \to \PP(X = 0)$, which is strictly less than $\sigma^{-K}$ by assumption \eqref{zeroatom}.  Let $r$ be large enough so that $M(-r) < \sigma^{-K}$.  Write $p = \sigma M(-r)^{1/K} < 1$, so that
			$$\PP(E_n) \le (p \E^{rA/B})^{Bn}.$$
		Choose the number $B$ to satisfy 
			\begin{equation} \label{defB}
				rA/\log(1/p) < B < rA/\log((1+p)/2p). \end{equation}
		This implies
			$$\tfrac{1+p}{2} < p \E^{r A/B} < 1.$$
		The left inequality will be used later in the proof of Thereom \ref{posconst}; the right inequality implies that
			$$\sum \PP(E_n) \le \sum (p \E^{r A/B})^{Bn} < \oo,$$
		so by the Borel-Cantelli lemma, with probability one, only finitely many of the events $E_n$ hold.
	\end{proof}

	\begin{proof}[Proof of Lemma \ref{upbound}]
		In what follows, we assume that $\Gamma$ is a $*$-connected set.  Fix an integer $B > 0$ and the sequence $a_n$.  Consider the events
			$$E_n = \{ \exists ~ \Gamma \mathrm{~such~that~} a_n \in \Gamma,~ X(\Gamma) > Cn,\mathrm{~and~} |\Gamma| \le Bn \}.$$
		We claim that we can choose $C$ large enough so that $\PP(E_n)$ decays exponentially in $n$.  From the Borel-Cantelli lemma, with probability one, it will follow that only finitely many of the events $E_n$ occur, which will prove the result.
	
		As in the proof of Lemma \ref{passstep}, we consider passage times of subsets of $\Gamma$ to exploit independence.  Unlike in that proof, it does not suffice to consider just one subset $\Gamma'$.  Instead, we partition the set into $k$ disjoint subsets $\Gamma_1, \dots, \Gamma_k$, and consider the passage times of each.
		
		There are $k = R^d$ points in the cube $\{0, \dots, R-1\}^d$; order them $z_1, \dots, z_k$.  We can partition the lattice $\Z^d$ into $k$ subsets by considering $R$-translations of these points.  We write $z = (z^1, \dots, z^d)$ for all $z \in \Z^d$.  For any $*$-connected set $\Gamma$, we partition it into $k$ subsets by defining
			$$\Gamma_j = \{ z \in \Gamma : \mbox{$R$ divides $z^i - z_j^i$ for all $i = 1, \dots, d$} \}$$
		for $j = 1, \dots, k$.  Since points in $\Gamma_j$ are separated by distance at least $R$, for each $j$ the family of random variables $\{X_z\}_{\Gamma_j}$ is independent.
		
		Let $\Gamma$ be a set described in the event $E_n$.  We may assume that $|\Gamma| = Bn$, since if it is less, the inclusion of additional points will only increase the passage time calculation.  Since
			$$X(\Gamma) = X(\Gamma_1) + \dots + X(\Gamma_k),$$
		if $X(\Gamma) > Cn$, then $X(\Gamma_j) > Cn/k$ for some $j$.  Furthermore, $|\Gamma_j| \le |\Gamma| \le Bn$ for each $j$.
		
		 As in Lemma \ref{passstep}, the probability $\PP(E_n)$ is bounded above by
			\begin{equation} \label{splitestimate2}
				\sum_{\Gamma} \PP \left( \exists ~ j \in \{1,\dots,k\} \mathrm{~s.t.~} X(\Gamma_j) > Cn/k \right), \end{equation}
		where again the sum is over $*$-connected sets $\Gamma$ for which $a_n \in \Gamma$ and $|\Gamma| = Bn$.  For each $j$, the family of random variables $\{X_i\}_{\Gamma_j}$ is independent.  Let $\{X_i\}$ be $Bn$ independent copies of $X$.  The exponential Chebyshev inequality \cite{durrett1988lnp} implies that the right-hand side of \eqref{splitestimate2} is bounded above by
			$$k \, \sigma^{Bn} \, \PP \left( \sum_{i=1}^{Bn} X_i > Cn/k \right) \le k \, \sigma^{Bn} \, \E^{-Cn/k} \, \EE \left( \E^{\sum X_i} \right),$$
		where we have added extra independent copies of $X_i$ to make the number of terms exactly $Bn$.  If we write $M = M(1) = \EE \E^X$, then this is equal to
			$$k \, \sigma^{Bn} \E^{-Cn/k} M^{Bn} = k \left( \sigma^{B} \E^{-C/k} M^{B} \right)^{n}.$$
		We choose $C \gg 1$ so that $\sigma^{B} \E^{-C/k} M^{B} < 1$.  Thus
			$$\sum \PP(E_n) \le k \sum \left( \sigma^{B} \E^{-C/k} M^{B} \right)^{n} < \oo,$$
		so by the Borel-Cantelli lemma, with probability one only finitely many of the events $E_n$ hold.
	\end{proof}

	\subsection{Applications to Continuum Model} \label{contapps}
		
	In this section, we apply the estimates from the previous section to the continuum model.  Lemma \ref{unifT} says that at large scales, the Riemannian distance function between two points is bounded by a uniform constant $K$ times the Euclidean distance between them.  Theorem \ref{posconst} says that $\mu_v$ is positive for all $v$.
	
	The cubes $C_z$ form a partition of $\R^d$.  For $x \in \R^d$, let $\tilde x$ be the unique point on the lattice such that $x \in C_{\tilde x}$.  The metric $g$ induces two notions of ``passage time'' over a discrete set $\Gamma$:
		$$\Lambda(\tilde\gamma) = \sum_{z \in \tilde\gamma} \Lambda_z \qquad \mathrm{and} \qquad \lambda(\tilde\gamma) = \sum_{z \in \tilde\gamma} \lambda_z.$$

	\begin{env_lem} \label{unifT}
		There exists a non-random $K > 0$ such that, with probability one, for all $x \in \R^d$ and $\rho > 0$, there exists $T(x) > 0$ such that if $t \ge T$ and $|x-y| \le \rho$, then
			$$d(tx,ty) \le Kt\rho.$$ 
	\end{env_lem}
	\begin{proof}
		By rescaling $x$ and $y$, it suffices to prove the lemma with $\rho = 1$.  
	
		Apply Lemma \ref{upbound} with $B = 2d$ and $X_z = \Lambda_z$.  Thus there exists a non-random $C > 0$ such that, with probability one, for any sequence $a_n \in \Z^d$, there exists $N > 0$ such that for all $n \ge N$, if $\Gamma$ is a finite $*$-connected set which contains the point $a_n$ and $|\Gamma| \le 2d n$, then $\Lambda(\Gamma) \le Cn$.
		
		Let $K = C\sqrt{d}$, and suppose that with positive probability, there is some $x \in \R^d$ such that for any $n > 0$, there exists $y$ (depending on $n$) such that $|x - y| \le 1$ but $d(nx, ny) > Kn$.  We will show that this leads to a contradiction.  Let $N$ be as in Lemma \ref{upbound} applied with the sequence $a_n = \widetilde{nx}$.
		
		Suppose that $n \ge N$.  Let $\gamma$ be the straight-line segment from $nx$ to $ny$.  Let
			$$\Gamma = \{ z \in \Z^d : \gamma \cap C_z \ne \emptyset \},$$
		index the cubes $C_z$ which $\gamma$ meets.  Note that $\widetilde{nx} \in \Gamma$.  Clearly, 
			$$|\Gamma| \le 2d |nx-ny| \le 2d n,$$
		since $|x-y| \le 1$.  By Lemma \ref{upbound},
			$$\Lambda(\Gamma) \le Cn.$$
		
		The distance between $nx$ and $ny$ is a lower bound for the Riemannian length of $\gamma$:
			$$d(nx, ny) \le R(\gamma).$$
		Furthermore, since $\gamma$ is a line segment, the Euclidean length of $\gamma$ in each cube is at most $\sqrt{d}$.  Thus we can estimate the Riemannian length of $\gamma$ by summing $\Lambda_z$ over $\Gamma$:
			\begin{equation} \label{partitionR}
				d(nx, ny) \le R(\gamma) = \sum_{z \in \Gamma} R(\gamma \cap C_z) \le \sum_{z \in \Gamma} \Lambda_z \sqrt{d} = \Lambda(\Gamma) \sqrt{d} \le C n \sqrt{d} = Kn, \end{equation}
		since $K = C\sqrt{d}$.  This contradicts the assumption that $d(nx,ny) > Kn$.
			\end{proof}

	\begin{env_thm} \label{posconst}
		The constants $\mu_v$ are all positive.
	\end{env_thm}
	\begin{proof}
		Suppose $\mu = \mu_v = 0$ for some unit vector $v$. 
		
		Let $\epsilon > 0$, and apply Lemma \ref{passstep} with the constant sequence $a_n \equiv 0$ to $A = 4\epsilon$ and $X_z = \lambda_z$.  Thus there exists a non-random $B > 0$ such that, with probability one, there exists $N_1 > 0$ such that for $n \ge N_1$, if $\Gamma$ is a $*$-connected set which contains the origin and $\lambda(\Gamma) \le 4\epsilon n$, then $|\Gamma| \le Bn$.  
		
		By Kingman's subadditive ergodic theorem, with probability one there exists $N_2 > 0$ such that if $n \ge N_2$, then
			\begin{equation} \label{kingmanposconst}
				d(0,nv) \le \epsilon n / 2, \end{equation}
		since we assumed that $\mu_v = 0$.
		
		Let $N = \max\{N_1, N_2\}$, and suppose $n \ge N$.  \emph{A priori}, the distance $d(0,nv)$ need not be realized as the Riemannian length of a curve, so let $\gamma$ be a $C^1$-curve from $0$ to $nv$ with
			$$R(\gamma) \le d(0,nv) + \epsilon n / 2 \le \epsilon n,$$
		where the second inequality follows from \eqref{kingmanposconst}.
	
		Define the discrete set
			\begin{equation} \label{Gammadef}
				\Gamma = \{ z \in \Z^d : L(\gamma \cap C_z) \ge 1/4 \}; \end{equation}
		that is, $z \in \Gamma$ provided the Euclidean length of $\gamma$ in the cube $C_z$ is at least $1/4$.  
		
		We claim that $\Gamma$ is $*$-connected.  Suppose not, so that the continuum set $W = \bigcup_{z \in \Gamma} C_z$ has at least two components, separated by Euclidean distance at least $1$.  Let $W'$ be the $1/4$-neighborhood around $W$, so that the components of $W'$ are separated by Euclidean distance at least $1/2$.  By definition of $\Gamma$, the curve $\gamma$ meets each component of $W'$, but not the complement $\R^d \smin W'$.  Since $\gamma$ is continuous, this is a contradiction; hence, $\Gamma$ is $*$-connected. 

		In each cube $C_z$, we can estimate the Riemannian length of $\gamma$ using $\lambda_z$:
			$$L(\gamma \cap C_z) \lambda_z \le R(\gamma \cap C_z),$$
		where $L$ denotes Euclidean length.  Furthermore, by summing $\lambda_z$ over the points of $\Gamma$, we get a lower bound for $R(\gamma)$:
			$$\tfrac{1}{4} \lambda(\Gamma) \le \sum_{z \in \Gamma} L(\gamma \cap C_z) \lambda_z \le \sum_{z \in \Gamma} R(\gamma \cap C_z) \le R(\gamma) \le \epsilon n.$$
		Clearly, $0 \in \Gamma$, so Lemma \ref{passstep} implies that $|\Gamma| \le Bn$.

		In the proof of that lemma, we chose $B$ so that
			$$B \le \tfrac{r}{\log((1+p)/2p)} A,$$
		for positive constants $r$ and $p < 1$ not depending on $A$; see \eqref{defB}.  Since $A = 4\epsilon$, if we write $B' = 4r/\log((1+p)/2p)$, then 
			\begin{equation} \label{GammaB}
				|\Gamma| \le B' \epsilon n. \end{equation}
			
		Let $z$ be a lattice point in $\Gamma$ which minimizes the distance $|nv - z|$.  Clearly, $|z| \ge n/2$.  Since $\Gamma$ is $*$-connected and $z \in \Gamma$,
			$$|\Gamma| \ge n/2\sqrt{d}.$$
		For small $\epsilon$, this contradicts \eqref{GammaB}, so $\mu$ must be positive.
	\end{proof}

	\section{The Shape Theorem and Consequences} \label{results}
	
	Define the function
		$$\mu(x) = \begin{cases} \mu_{x/|x|} |x|, & x \ne 0 \\ 0, & x = 0. \end{cases}$$
	Proposition \ref{contmu} and Theorem \ref{posconst} imply that $\mu$ is continuous and, for $x \ne 0$, strictly positive.  It follows from the triangle inequality for the distance function that $\mu(x)$ is a norm on $\R^d$.  Consider the unit ball in this norm,
		\begin{equation} \label{limshape}
			A = \{ x : \mu(x) \le 1 \} = \{ x : |x| \le \mu_{x/|x|}^{-1} \}, \end{equation}
	as well as the random Riemannian ball of radius $t$ centered at the origin,
		$$B_t = \{ x : d(0,x) \le t \}.$$
		
	\begin{env_thm}[Shape Theorem]	\label{shapethm}
			For all $\epsilon > 0$, with probability one, there exists $T$ such that if $t \ge T$, then
				\begin{equation} \label{shapestatement}
					(1-\epsilon) A \subseteq \tfrac{1}{t} B_t \subseteq (1+\epsilon) A.
				\end{equation}  
	\end{env_thm}
	
	The set $A$ is called the limiting shape of the random Riemannian metric $g$.  The shape theorem for lattice first-passage percolation was proved by Cox and Durrett \cite{cox1981slt}; our proof is modeled on the arguments in Durrett \cite{durrett1988lnp}.
		
		\begin{proof}
			It suffices to prove the theorem for $\epsilon \in (0,1)$.  Let
				\begin{equation} \label{delta}
					\delta < \min \left\{ \frac{1}{K}, 1 - \frac{1+\epsilon^2}{1+\epsilon} \right\}, \end{equation}
			where $K$ is as in Lemma \ref{unifT}.  Let $\Beuc(x,r)$ denote the Euclidean ball of radius $r$ centered at $x$.
			
			We will show that with probability one, there exists $T > 0$ such that if $t \ge T$, then $(1-\epsilon)A \subseteq \tfrac{1}{t} B_t$.  To do this, we will first prove that for every $x$, there is a random $T(x)$ such that if $t \ge T$, then the small Euclidean ball $\Beuc(x,\delta\epsilon^2)$ is contained in $\tfrac{1}{t} B_t$.  Since the set $(1-\epsilon)A$ is compact, it can be covered by finitely many balls $\Beuc(x_i, \delta\epsilon^2)$.  Letting $T = \max\{T(x_i)\}$ will prove the result.
			
		Fix $x \in (1-\epsilon)A$.  We claim that with probability one, there exists $T(x)$ such that if $t \ge T$ and $|x-y| \le \delta \epsilon^2$, then $d(0,ty) \le t$, hence $\Beuc(x,\delta \epsilon) \subseteq \tfrac{1}{t} B_t$ for $t \ge T$.  By the triangle inequality,
			$$d(0,ty) \le d(0,tx) + d(tx,ty).$$
		
		The first term is controlled by Kingman's theorem:  with probability one, there exists $T_1(x)$ such that if $t \ge T_1$, then
			$$d(0,tx) \le (1+\epsilon) t \mu(x) \le (1-\epsilon^2) t,$$
		since $\mu(x) \le (1-\epsilon)$.
	
		The second term is controlled by Lemma \ref{unifT} applied to this $x$ and $\rho = \delta\epsilon^2$.  With probability one, there exists $T_2(x)$ such that for all $t \ge T_2$ and $y$ with $|x-y| \le \delta \epsilon^2$, then
			$$d(tx,ty) \le Kt \delta \epsilon^2 < \epsilon^2 t,$$
		since $\delta < 1/K$.
		
		Let $T(x) = \max\{T_1, T_2\}$.  If $t \ge T(x)$, then for all $y \in \Beuc(x,\delta \epsilon^2)$, $$d(0,ty) \le (1-\epsilon^2)t + \epsilon^2 t = t,$$ implying that $\Beuc(x, \delta \epsilon^2) \subseteq \tfrac{1}{t} B_t$ for all $t \ge T(x)$.  Let $x_i$ be finitely many points such that $\bigcup \Beuc(x_i, \delta \epsilon^2)$ covers $(1-\epsilon)A$, and let $T = \max_i T(x_i)$.  Then for all $t \ge T$,
			$$(1-\epsilon)A \subseteq \bigcup_i \Beuc(x_i, \delta \epsilon^2) \subseteq \tfrac{1}{t} B_t,$$
		completing the lower half of the shape theorem.
		
		Now we prove the upper half.  For any $x$, let $T_2(x)$ be as above, so that if $t \ge T_2$ and $|x-y| < \delta \epsilon^2$, then $d(tx,ty) \le \epsilon^2 t$.  By Kingman's theorem, with probability one, there exists  $T_3(x)$ such that if $t \ge T_3$, then
			$$d(0,tx) \ge  (1-\delta) t \mu(x).$$
		Choose finitely many $x_i \in 2A \smin (1+\epsilon)A$ such that the closure of $2A \smin (1+\epsilon)A$ is covered by $\bigcup \Beuc(x_i, \delta \epsilon^2)$.  Belonging to $2A \smin (1+\epsilon)A$ implies that $\mu(x_i) > (1+\epsilon)$.  Let $T = \max_i \{T_2(x_i), T_3(x_i)\}$, and let $t \ge T$.  Then if $y \in \Beuc(x_i, \delta \epsilon^2)$,
			\begin{eqnarray*}
			 	d(0,ty) &\ge& d(0, tx_i) - d(tx_i, ty) \\
					&\ge& (1-\delta)t \mu(x_i) - \epsilon^2 t \\
					&\ge& (1-\delta)(1+\epsilon)t - \epsilon^2 t \\
					&>& (1+\epsilon^2) t - \epsilon^2 t = t,
			\end{eqnarray*}
		where the final inequality follows from the assumption \eqref{delta} on $\delta$.  Thus $\bigcup \Beuc(x_i, \delta \epsilon^2) \subseteq \tfrac{1}{t} B_t^c$ so
			$$2A \smin (1+\epsilon)A \subseteq \bigcup \Beuc(x_i, \delta \epsilon^2) \subseteq \tfrac{1}{t} B_t^c.$$
		The set $\tfrac{1}{t} B_t$ is connected and contains the origin, hence $\tfrac{1}{t} B_t \subseteq (1+\epsilon)A$ as desired.
	\end{proof}
	
	\begin{env_cor} \label{isotropic}
		If the measure $\PP$ is isotropic (rotationally-invariant), then the limiting shape $A$ is the Euclidean ball of radius $\mu^{-1}$. \hfill $\square$
	\end{env_cor}
	
	A simple consequence of Lemma \ref{unifT} is that the convergence \eqref{kingman} given by Kingman's theorem is uniform:
	
	\begin{env_pro} \label{unifconv}
%		For all $\epsilon > 0$, there exists an $M$ such that if $|x| \ge M$, then
%			$$|d(0,x) - \mu(x)| \le \epsilon |x|.$$
		For all $\epsilon > 0$, with probability one, there exists $T>0$ such that if $t \ge T$, then for all $v \in S^{d-1}$,
			$$\left| \tfrac{1}{t} d(0,tv) - \mu_v \right| \le \epsilon.$$
	\end{env_pro}
	
	\begin{proof}
		Suppose not.  Then with positive probability, there exist $\epsilon > 0$, $t_n \to \oo$ and $v_n \in S^{d-1}$ such that
			$$\left| \tfrac{1}{t_n} d(0,t_n v_n) - \mu_{v_n} \right| > \epsilon,$$
		for all $n$.  By compactness of the sphere, a subsequence of $v_n$ converges to some $v \in S^{d-1}$; assume without loss of generality that $v_n \to v$.  By the triangle inequality,
			\begin{equation} \label{epsunif}
				\epsilon < \left| \tfrac{1}{t_n} d(0,t_nv_n) - \tfrac{1}{t_n} d(0, t_nv) \right| + \left| \tfrac{1}{t_n} d(0,t_nv) - \mu_v \right| + \left| \mu_{v_n} - \mu_v \right|. \end{equation}
		The second and third terms tend to zero a.s. as $n \to \oo$ by Kingman's theorem and Proposition \ref{contmu}, respectively.  By the triangle inequality, the first term is bounded by $\tfrac{1}{t_n} d(t_nv, t_nv_n)$.  Let $\rho = \epsilon / 2K$, where $K$ is as in Lemma \ref{unifT}.  Since for large $n$, $|v - v_n| \le \rho$, we can apply that lemma with $x = v$ to get
			$$\tfrac{1}{t_n} d(t_nv,t_nv_n) \le K\rho = \epsilon /2,$$
		for large $n$ almost surely.  This contradicts \eqref{epsunif}.
	\end{proof}

	\begin{env_thm} \label{completeness}
	 	With probability one, if $\gamma$ is a $C^1$-curve parametrized by Riemannian length, then $|\gamma(t)| < \oo$ for all $t \ge 0$. 
	\end{env_thm}
	\begin{proof}
		It suffices to consider curves starting from the origin.  Suppose that with positive probability, there exists a smooth curve $\gamma$ parametrized by Riemannian length which starts from the origin and for which
			\begin{equation} \label{Tlim}
				\lim_{t \to T} |\gamma(t)| = \oo \end{equation}
		for some finite $T$.  Since $\gamma$ is parametrized by Riemannian length, for all $t \ge 0$,
			$$t \ge d(0,\gamma(t)).$$
		Proposition \ref{unifconv} and \eqref{Tlim} imply that
			$$0 = \lim_{t\to T} \frac{t}{|\gamma(t)|} \ge \lim_{t\to T} \frac{d(0,\gamma(t))}{|\gamma(t)|} \ge \min_{v} \mu_v > 0,$$
		a contradiction.
	\end{proof}
	
	We can refine this result if we impose an additional smoothness constraint on the metric.  If $g$ is a $C^2$-smooth Riemannian metric, i.e. $g \in C^2(\R^d, \SPD)$, then we can use the calculus of variations to derive the Euler-Lagrange equations for the functional $R$.  These are called the geodesic equations \cite{lee1997rmi} for the Riemannian metric $g$, and solutions to this system are called geodesics.  The geodesic equations form a second-order system with locally-Lipschitz coefficients, so a geodesic is uniquely determined by its starting point and velocity.  We call a geodesic $\gamma$ length-minimizing if for all $x, y \in \gamma$, the distance $d(x,y)$ is realized as the Riemannian length of the part of the curve $\gamma$ which connects the two points.  Not all geodesics are length-minimizing; for example, on the sphere, the geodesics are great circles, which do not minimize length past antipodal points.

	A metric is said to be geodesically complete if for all $x \in \R^d$ and $v \in T_x \R^d$, the unique geodesic $\gamma$ at $x$ in direction $v$ can be continued for all time.  Part of the Hopf-Rinow theorem \cite{lee1997rmi} of Riemannian geometry is that geodesic completeness is equivalent to the condition stated in Theorem \ref{completeness} for our random metric $g$.  A further corollary \cite{lee1997rmi} is that distances are always realized by geodesics.  Summarizing, we have
	
	\begin{env_cor} \label{realized}
		Suppose that, in addition to the assumptions of this paper, $g$ is a $C^2$-smooth random Riemannian metric.  Then, with probability one, $g$ is geodesically complete.  Consequently, for all $x, y \in \R^d$, there is a finite length-minimizing geodesic $\gamma$ connecting $x$ to $y$ such that
			$$d(x,y) = R(\gamma).$$
	\end{env_cor}

We alert the reader to a different meaning of the word ``geodesic,'' used often in the first-passage percolation literature.  The term is used there to denote a globally length-minimizing path.  This is very different from the standard meaning of the word in differential geometry: as described above, geodesics are the curves which locally minimize length, but not necessarily globally.  We adhere to this meaning in the present paper.

	The existence of two-sided minimizing paths is an open question for all FPP models.  For two-dimensional standard FPP, Licea, Newman and Piza \cite{licea1996geodesics, licea1996superdiffusivity, newman1997tds} have a number of results in this direction.  Their work relies on certain curvature assumptions about the limiting shape, which have not been verified for models of independent FPP.  Howard and Newman's \cite{howard1997euclidean, howard2001special} model of Euclidean FPP, on the other hand, is rotationally invariant.  Consequently, the limiting shape is a Euclidean ball, and they prove many results not available in the lattice setting.  See the excellent survey \cite{howard2004mfp} for more details.  By Corollary \ref{isotropic}, the limiting shape for isotropic Riemannian FPP is a Euclidean ball.  Our hope is for this setting to be a fertile ground for adapting the results referenced in this paragraph.\newline

	{\bf Remark added after publication:}  The moment assumptions \eqref{finmom} and \eqref{finmom2} can be slightly weakened to:  $$M(r) = \EE[ \E^{r \Lambda_0} ] < \oo \qquad \mbox{for all $r \le a$},$$ for some $a > 0$.  The only modification to the paper is in the proof of Lemma \ref{passstep}, where one uses $M = M(a)$ instead of $M = M(1)$.  We leave the details to the reader.
	
%\bibliographystyle{plain}
%\bibliography{biblio}
\end{document}